\documentclass[bj,preprint]{imsart}

\setattribute{journal}{name}{}
\usepackage[utf8]{inputenc}
\usepackage[T1]{fontenc}
\usepackage[english]{babel}
\usepackage{graphicx}
\usepackage{xspace}
\usepackage{lmodern}
\usepackage{amsmath,amssymb,amsfonts,amsthm,shuffle}
\usepackage{enumerate}
\usepackage{tikz}\usetikzlibrary{arrows}

\RequirePackage[OT1]{fontenc}
\RequirePackage{amsthm,amsmath}
\RequirePackage[numbers]{natbib}
\RequirePackage[colorlinks,citecolor=blue,urlcolor=blue]{hyperref}

\DeclareMathOperator{\Var}{Var}

\newcommand{\prob}[1]{\mathbb{P}\left(#1\right)}

\newcommand{\indic}[1]{\mathrm{\textbf{1}}_{#1}}

\newcommand{\Esp}[1]{\mathbb{E}\left[#1\right]}
\newcommand{\Espc}[2]{\mathbb{E}\left[ #1 \middle\vert #2 \right]}
\newcommand{\Espi}[2]{\mathbb{E}_{#1}\left[ #2  \right]}
\newcommand{\Espci}[3]{\mathbb{E}_{#1}\left[ #2 \middle\vert #3 \right]}

\newcommand{\probc}[2]{\mathbb{P}\left( #1 \middle\vert #2 \right)}

\newcommand{\setC}{\mathbb{C}}
\newcommand{\setN}{\mathbb{N}}
\newcommand{\setZ}{\mathbb{Z}}

\newcommand{\setR}{\mathbb{R}}

\newcommand{\ca}[1]{\mathcal{#1}}

\newcommand{\ha}[1]{\widehat{#1}}

\DeclareMathOperator{\card}{card}

\newcommand{\norm}[1]{ \left\vert\left\vert #1 \right\vert\right\vert }

\newcommand{\IntP}[1]{\lfloor #1 \rfloor}
\newcommand{\FracP}[1]{ \{#1\}}

\newcommand{\qshuffle}{\ha{\shuffle}}

\newtheorem*{theorem*}{Theorem}
\newtheorem{theorem}{Theorem}[section]

\newtheorem{definition}{Definition}[section]

\newtheorem{property}{Property}[section]
\newtheorem{lemma}{Lemma}[section]
\newtheorem{proposition}{Proposition}[section]

\begin{document}

\begin{frontmatter}
\title{Area anomaly in the rough path Brownian scaling limit of hidden Markov walks}
\runtitle{Area anomaly for hidden Markov walks}

\begin{aug}
\author{\fnms{Olga} \snm{Lopusanschi}
\ead[label=e1]{olga.lopusanschi@gmail.com}}
\and
\author{\fnms{Damien} \snm{Simon}\ead[label=e2]{damien.simon@sorbonne-universite.fr}}

\runauthor{O. Lopusanschi and D. Simon}

\affiliation{Laboratoire de Probabilit\'es, Statistique et Mod\'elisation (UMR~8001), Sorbonne Universit\'e, CNRS, France.}

\address{Laboratoire de Probabilit\'es, Statistique et Mod\'elisation (UMR~8001),\\ 
Sorbonne Universit\'e and CNRS,\\ 
4 place Jussieu, 75005 Paris, France.\\
\printead{e1},\\
\printead*{e2}
}
\end{aug}


\begin{abstract}
We study the convergence in rough path topology of a certain class of discrete processes, the hidden Markov walks, to a Brownian motion with an area anomaly. This area anomaly,
which is a new object, keeps track of the time-correlation of the discrete models and brings into light the question of embeddings of discrete processes into continuous time. We also identify an underlying combinatorial structure in the hidden Markov walks, which turns out to be a generalization of the occupation time from the classical ergodic theorem in the spirit of rough paths. 
\end{abstract} 

\end{frontmatter}

\tableofcontents

\section{Introduction}
\subsection{The context}

Rough paths theory was introduced by T.~Lyons in 1998 (see, for example, \citep{LyonsOriginal,StFlour}) in order to provide a deterministic setting to stochastic differential equations (SDEs) of the type \[
dy_t = f(y_t)[dx_t]
\]
where $(y_t)_t$ is a path in a finite-dimensional vector space $V'$, $(x_t)_t$ is a path in another finite-dimensional vector space $V$ with H\"older regularity $\alpha<1$ (which is often the case for stochastic processes) and $f: V' \to Hom(V,V')$ is a smooth map. Whenever the classical Young integration \citep{bookYoung} fails (which is the case for $\alpha<1/2$), paths may be lifted (in a non-unique way) to a larger more abstract space, the \emph{rough path space}, for which existence, uniqueness and continuity of the solution map hold and become easier to prove.

The rough path space of level $k\geq 1$ corresponds to paths with values in \[
T_1^{(k)}(V)=V\oplus V^{\otimes 2} \oplus \ldots V^{\otimes k}
\]
and a finite-variation path $x:[0,T]\to V$ is lifted to a path $S(\gamma):[0,T]\to T_1^{(k)}(V)$, called the step-$k$ signature of the path, through the formulae:
\begin{align*}
S(x)(t)&= \left( S_1(x)(t),\ldots,S_k(x)(t)\right) 
\\
S_j(x)(t) &= \int_{0<s_1<\ldots <s_j < t} dx_{s_1}\otimes dx_{s_2}\otimes \ldots \otimes dx_{s_j}
\end{align*}
The coefficients of the signature satisfy a set of algebraic relations (shuffle product and concatenation product) and of analytic bounds (depending on the H\"older regularity $1/2<\alpha$ of $x$), which will be detailed below. For a regularity $\alpha<1/2$, some of the integrals above are not well-defined any more but the algebraic relations and the analytic bounds are kept as a new \emph{definition} of the signature, and therefore also defines formally the iterated integrals.

\paragraph{Some examples of rough paths.} The space of rough paths of level 2, for example, contains the enhanced Brownian motion \citep{FrizVicBook,Breuillardetal} \begin{equation}\label{eq:enhancedBM}
t\mapsto (B_t, (1/2)B_t\otimes B_t+ A_t^{\text{L\'evy}})
\end{equation} where $(B_t)_t$ is a $V$-valued Brownian motion and $A_t^{\text{L\'evy}}$ its L\'evy area (with values in $V\wedge V$). The Brownian motion has regularity $1/2^-$ and its L\'evy area is defined through the usual stochastic calculus by
\[
A_t^{\text{L\'evy}} = \frac{1}{2} \int_{0<s_1<s_2<t} (dB_{s_1} \otimes dB_{s_2} - dB_{s_2}\otimes dB_{s_1})
\] 
where the stochastic integration may be either in the It\^o or  Stratonovich sense.

The space of rough paths also contains less trivial objects such as two-dimensional area bubbles \citep{artLejayLyonsArea}, defined as the limit signature of the sequence of paths \[
x_n(t) = \frac{1}{\sqrt{n}} \left( \cos(nt), \sin(nt) \right)
\]
which turns around $(0,0)$ faster as $n$ increases. The signature is given by $((0,0),(t/2)e_1\wedge e_2)$ (where $(e_1,e_2)$ is the canonical basis of $\setR^2$): the first level is constant at $(0,0)$, which corresponds to a constant path at $(0,0)$, whereas the second level corresponds to an accumulation of area at constant speed. This shows that the rough path space contain more information than the usual path space about the microscopic structure of approximations.

A quick look at the rough path space also shows that a combination of the two previous examples, i.e. the signature
\begin{equation}\label{eq:enhancedanomalousBM}
S(B)(t)=\left( B_t, (1/2)B_t\otimes B_t+ A_t^{\text{L\'evy}} + t \Gamma \right),
\end{equation}
where $(B_t)$ is a two-dimensional Brownian motion and $\Gamma$ is a fixed element of $V\wedge V$,
is a valid rough path. From the point of view of stochastic integration however, this object is new since the construction of the L\'evy area does not depend on the choice of It\^o or Stratonovich integration and thus the $\Gamma$ term is not an artefact of the choice of the stochastic integral. Thus, one may wonder whether this type of objects is relevant for applications and  continuous time limits of discrete models.

\paragraph{Continuous time limit of discrete models in the rough path space: existing results.} A first result is the generalization of Donsker theorem for i.i.d. random variables by Breuillard et al. \citep{Breuillardetal}, which, without surprise, gives the convergence in rough path topology to the enhanced Brownian motion \eqref{eq:enhancedBM} under suitable finite moment hypotheses.

A first result of convergence to \eqref{eq:enhancedanomalousBM} with non-zero $\Gamma$ was obtained in \citep{LopusanschiSimon}. It relies mostly on a geometric construction of periodic graphs and a non-reversible random walk dynamics on the graph. The idea is the following: at the scale of the period of the graph, the dynamics looks like the one of a random walk but, at the smallest scale of the periodic pattern, windings may occur and contribute to the area anomaly $\Gamma$. 

\paragraph{The generalization to hidden Markov walks.} The present paper presents a much more general construction of discrete time models converging to \eqref{eq:enhancedanomalousBM} in the rough path topology. In particular, the construction does not rely anymore on geometric properties, is purely probabilistic and the emergence of a non-zero $\Gamma$ relies only on probabilistic short time correlations of the discrete time model. 

Many discrete models may enter the present settings and exhibit a non-zero area anomaly in the limit. It was indeed a surprise to us that no such result have been previously considered in the literature, despite the vast literature about stochastic integration, continuous-time limits and rough paths.

The motivation for considering hidden Markov walks as defined below is multiple. When considering limits such as \eqref{eq:enhancedanomalousBM},  it is very tempting to consider ad hoc models obtained by a discretization of the limit but it does not shed any new light on the model and one of our purposes was to avoid such an approach. The central limit theorem and its process extension, Donsker theorem, show that obtaining anomalous behaviour around a normal law/Brownian motion from Markov models requires a fine tuning of the correlations at the discrete level. This is precisely what the present paper describes through hidden Markov walk. 

Moreover, hidden Markov walks are simple objects from the point of view of numerical simulations and thus the present construction may lead to interesting discrete approximations, as for example in \citep{artDavieDiscrete,artDiscreteIntegrals}.

Finally, the proofs of the results presented below introduce new objects in the study of discrete time Markov chains: iterated occupation times, which are a discrete analogue of iterated integrals and signature of the rough path theory. These iterated occupation times satisfy algebraic relations, such as quasi-shuffle and concatenation, and are of independent interest in the general theory of Markov chains.

\subsection{Formulation of the results}

\subsubsection{Hidden Markov walks and the first theorem.}

We first define the discrete time model we will consider throughout the paper and we choose to call it \emph{hidden Markov walk}, which is a particular case of hidden Markov chain as introduced in \citep{firstHMC}.

\begin{definition}[hidden Markov walk]\label{def:hiddenmarkovwalk} Let $E$ be a countable set and $V$ a finite dimensional real vector space. A \emph{hidden Markov walk} is a process $(R_n,X_n)_{n\in\setN}$ on $E\times V$ such that: 
\begin{enumerate}
\item the process $(R_n)_{n\in\setN}$ is a Markov chain on $E$;
\item conditionally on the process $(R_n)$, the increments $F_n=X_{n+1}-X_{n}$ are independent and have marginal laws such that, for any Borel set $A$ of $V$, any $n\in\setN$ and any $r\in E$, it holds:
\begin{align}
\probc{ F_n\in A}{\sigma(R)} &= \probc{ F_n\in A}{R_n}, 
\\
\probc{ F_n\in A}{R_n=r} &= \probc{ F_1\in A}{R_1=r}.
\end{align}
\end{enumerate}
\end{definition}

The process $(R_n,F_n)_{n\in\setN}$ corresponds to the usual definition of a hidden Markov chain. The additional vector space structure of $V$ allows one to consider the $(F_n)$ as increments and to add them to obtain the process $(X_n)$. 

A classical way of embedding the discrete process $(X_n)_{n\in\setN}$ in continuous time is the so-called Donsker embedding.
\begin{definition}[Donsker embedding]
Let $N\in\setN$ and $V$ be a vector space. Let $(x_n)_{0\leq n\leq N}$ be a $V$-valued sequence. The Donsker embedding $\iota_N(x)$ is the path $\iota_N(x):[0,1]\to V$ such that, for any $0\leq k< N$ and any $s\in [0,1]$, 
\[
\iota_N(x)\left( (1-s)\frac{k}{N}+s\frac{k+1}{N} \right)
=
(1-s) x_k + s x_{k+1},
\]
which interpolates linearly between the $(x_n)$ and accelerate time by a factor $N$.
\end{definition}

We also define, for any real number $s>0$, the dilation operators:
\begin{align*}
\delta_s : T_1^{(n)}(V) & \to T_1^{(n)}(V)
\\
u_1\otimes\ldots \otimes u_l  & \mapsto s^l (u_1\otimes\ldots \otimes u_l)
\end{align*}
where the $(u_i)_{1\leq i\leq l}$ are vectors of $V$ and $l$ an integer smaller than $n$.

We may now state our first theorem of convergence to an anomalous enhanced Brownian motion. The precise definition of the topological rough paths space $\ca{C}^{\alpha}([0,1],G^2(V))$ is given in section \ref{sec:introroughpath}.

\begin{theorem}\label{thm:convergenceDonsker}
Let $(R_n,X_n)_{n\in\setN}$ be a hidden Markov walk on $E\times V$ such that, for a fixed $r_0\in E$,
\begin{enumerate}[(i)]
\item $X_0=0$ a.s. and $R_0=r_0$ a.s.
\item the Markov chain $(R_n)$ is irreducible and positive recurrent, with invariant probability $\nu$, and the first return time
\begin{equation}
T_1=\inf\{n\geq 1; R_n=R_0 \}
\end{equation}
has finite moments of all orders.
\item the increments $F_n=X_{n+1}-X_n$ satisfy, for all $p\in\setN$,
\begin{equation}
\sup_{r\in E} \Espc{ \norm{F_1}^p}{R_1=r} < \infty
\end{equation}
\item the walk $(X_n)$ is centred, i.e. satisfies
\begin{equation}
\sum_{r\in E} \nu(r) \Espc{F_1}{R_1=r} = 0_V 
\end{equation}
\end{enumerate}
Then, for any $1/3<\alpha<1/2$,
the sequence of $G^2(V)$-valued continuous time processes $(\delta_{N^{-1/2}}\circ S\circ \iota_N(X)(t))_{t\in[0,1]}$ (where $S$ is the signature of the path) converges in law in the rough path topology of $C^{0,\alpha\text{-Höl}}([0,1],G^2(V))$ to the enhanced Brownian motion (as defined in~\eqref{eq:enhancedanomalousBM}) with covariance matrix $C$ and anomalous area drift $\Gamma$. Moreover, the limit law does not depend on the choice of the initial law of $R_0$.

Explicit formulae for $C$ and $\Gamma$ are presented in \eqref{eq:formulacovariancematrixC}, \eqref{eq:formulaareadrift}, \eqref{eq:formulaforCocctimes} and \eqref{eq:formulaforGammaocctimes}.
\end{theorem}

The enhanced Brownian motion obtained in the limit is not the standard one, with identity covariance matrix and zero area drift, but has a covariance $\Esp{B_t^{(i)} B_t^{(j)}}=C_{ij} t$ and an area $\Esp{A^{(ij)}_t}= t\Gamma_{ij}$ (the standard L\'evy area is centred). As detailed during the proof of the theorem and in section~\ref{sec:studyCGamma}, the main ingredient for the existence of a non-zero $\Gamma_{ij}$ is the fact that the area covered during one excursion of the Markov chain $(R_n)$ may be \emph{not} centered and we have:
\begin{equation}
\Espi{r_0}{ A^{(ij)}_{T_1} }  = \Gamma_{ij} \Espi{r_0}{T_1}
\end{equation}
where $T_1$ is the length of the excursion ans $A^{(ij)}_{T_1}$ is the area covered by the discrete walk $X_n$ between $0$ and $T_1$.

The first hypothesis is not restrictive at all since there is global translation invariance. The irreducibility of $(R_n)$ is not restrictive since one may always restrict $E$ to one of its irreducible component. The hypotheses of positive recurrence and of finite moments of $T_1$ are important for the proofs but are trivially satisfied whenever $E$ is finite. The hypothesis on the moments of the increments is already required in \citep{Breuillardetal} to obtain rough path convergence of random walks. The assumption on centring is not restrictive: the centring is important to describe TCL-like fluctuations around the law-of-large-number asymptotics of $X_n$.

\subsubsection{Iterated occupation times.}

\paragraph{Definition and basic algebraic properties.}
The proof of theorem \eqref{thm:convergenceDonsker} uses extensively conditional expectations with respect to the process $\sigma(R)$ and the hidden Markov structure reduces iterated integrals with values in tensor products of the space $V$ to combinatorial quantities related to the single process $(R_n)_{n\in\setN}$. It appears that these combinatorial quantities also have a nice algebraic structure and interesting asymptotics.

\begin{definition}[iterated occupation time of a sequence]
Let $E$ be a set and let $(x_n)_{0\leq n<N}$ be a sequence of elements of $E$. For any $k\geq 1$ and any sequence $u=(u_1,\ldots,u_k)\in E^k$, the iterated occupation time of $(x_n)$ at $u$ is defined as:
\begin{equation}
\label{eq:defiteratedocctime}
\begin{split}
L_{u}(x) = 
\card\{
& (i_1,\ldots,i_k)\in \{0,...,N-1\}^k ; 
\\
&i_1<i_2<\ldots<i_k \quad\text{and}\quad x_{i_1}=u_1, x_{i_2}=u_2,\ldots,x_{i_k}=u_k 
\}
\end{split}
\end{equation}
and, if $k=0$, the empty sequence is noted $\epsilon$ and by convention $L_\epsilon(x)=1$.
\end{definition}

For $k=1$, $L_u(x)$ counts the number of times the sequence $x$ visits $u$, hence the name of iterated occupation times for larger $k$. Another writing of \eqref{eq:defiteratedocctime} makes the relation with iterated integrals clearer:
\begin{equation}
L_{u}(x)
= \sum_{0\leq i_1<\ldots <i_k <N}
\indic{x_{i_1}=u_1}\ldots\indic{x_{i_k}=u_k}
\end{equation}
The iterated sum structure endows the iterated occupation times with both a concatenation structure and a quasi-shuffle structure.
\begin{property}[concatenation structure]
Let $(x_n)_{0\leq n<N}$ be a finite sequence of elements of $E$ of length $N$. Let $u=(u_1,\ldots,u_k)\in E^k$ for some $k$. Let $M$ be an integer smaller than $N$. Then it holds
\begin{equation}
L_{u_1\ldots u_k}(x) = \sum_{p=0}^k L_{u_1\ldots u_p}( (x_n)_{0\leq n<M}) L_{u_{p+1}\ldots u_k}( (x_n)_{M\leq n< N})
\end{equation}
where, by convention, $L_\epsilon(x)=1$ for $\epsilon$ the zero length sequence.
\end{property}

The quasi-shuffle property requires some additional combinatorial definitions. There are various approaches to quasi-shuffles and the interested reader may refer to \citep{artQuasiShuffle,artQuasiHoff2} for more results on this notion.
\begin{definition}[quasi-shuffle product]
Let $E$ be a set and let $A(E)$ be the algebraic direct sum \[
A(E) = \oplus_{k\geq 0} \oplus_{u\in E^k} \setR u
\]
(for $k=0$, the empty sequence is written $\epsilon$) endowed with the shift operators $a\cdot u$ defined for $a\in E$ and $u=(u_1,\ldots,u_k)\in E^k$ by
\begin{align*}
a\cdot u &= (a,u_1,\ldots,u_k)
\end{align*}
The quasi-shuffle product $\qshuffle$ is defined recursively on the canonical basis by:
\begin{align*}
\epsilon \qshuffle \epsilon &= \epsilon \\
(a\cdot u) \qshuffle \epsilon &= a\cdot u \\
\epsilon \qshuffle (b\cdot v)  &= b\cdot v \\
(a\cdot u) \qshuffle (b\cdot v) 
&= a\cdot (u\qshuffle (b \cdot v)) + b \cdot ((a\cdot u)\qshuffle v) \\
&\phantom{=} + \indic{a=b} a\cdot (u\qshuffle v)
\end{align*}
and extended to $A(E)$ by linearity.

To any finite sum $a=\sum_{p} c_p u_p$ with elements $u_p\in \cup_k E^k$, the definition of the iterated occupation time $L_u(x)$ of a sequence is extended by:
\begin{equation}\label{eq:extensionquasishuffle}
L_a(x) =\sum_{p} c_p L_{u_p}(x)
\end{equation}
and is thus a linear application $A(E)\to\setR$.
\end{definition}
When the sequences $u_1$ and $u_2$ have no element in common, the set of quasi-shuffles $u_1\qshuffle u_2$ is equal to the classical shuffle of the two sequences. We now express the fact that the linear map $L_\bullet : A(E) \to\setR$ is a morphism of algebra for the quasi-shuffle product.

\begin{proposition}\label{prop:occtime:quasishuffle} Let $(x_n)_{0\leq n<N}$ be a finite sequence of elements of $E$. For any $k,l\geq 1$ and any sequences $u=(u_1,\ldots,u_k)\in E^k$ and $v=(v_1,\ldots,v_l)\in E^l$, it holds:
\begin{equation}
L_u(x) L_v(x) = L_{u\qshuffle v}(x)
\end{equation}
where the r.h.s has to interpreted in $A(E)$ through the extension~\eqref{eq:extensionquasishuffle}.
\end{proposition}
Checking the following identity for $k=l=1$ is left as a warm-up exercise for the complete proof in section~\ref{sec:occtimeproof}
\begin{align*}
L_{(u_1)}(x) L_{(v_1)}(x)
= L_{(u_1v_1)}(x) + L_{(v_1u_1)}(x) +
\begin{cases}
L_{(u_1)(x)} & \text{if $u_1=v_1$} \\
0 & \text{else}
\end{cases}
\end{align*}
The replacement of shuffles of iterated integrals by quasi-shuffles of iterated sums is related to the fact that the Lebesgue measure puts zero mass on singlets nor on $d-1$-dimensional subspaces of $\setR^d$. 

\paragraph{From hidden Markov paths to iterated occupation times.}
As announced, the definition~\ref{def:hiddenmarkovwalk} of hidden Markov walks fits nicely with the framework of iterated occupation times.

\begin{property}
Let $(R_n,X_n)$ be a hidden Markov walk on $E\times V$. Let $N\in\setN$ and $\mathbf{X}_N : [0,1]\to V$ be the piecewise linear path $\mathbf{X}_N=\iota_N(X)$. Then, under suitable integrability hypothesis for the existence of the conditional expectation, there exist deterministic coefficients $(f_u)_{u\in E}$ in $V$, $(c_u)_{u\in E}$ in $V\otimes V$ and $(b_{u,v})_{u,v\in E}$ in $V\otimes V$ such that
\begin{align*}
\Espc{ \int_0^1 d\mathbf{X}_N(s)  }{ \sigma(R) } &= \sum_{u\in E} f_u L_u( (R_n)_{0\leq n<N})
\\
\Espc{ \int_{0<s_1<s_2<1} d\mathbf{X}_N(s_1) \otimes d\mathbf{X}_N(s_2)  }{ \sigma(R) } &= \sum_{u\in E} c_u L_u( (R_n)_{0\leq n<N})
\\
& \phantom{=} + \sum_{(u,v)\in E} b_{u,v} L_{uv}( (R_n)_{0\leq n<N})
\end{align*}
and more generally such decompositions on iterated occupation time hold for higher iterated integrals.
\end{property}
One checks that iterated integrals of level two are mapped to iterated occupation times of both level one and level two: the absence of coherent grading is related to the fact that iterated integrals satisfy a shuffle  property whereas iterated sums satisfy only a quasi-shuffle property.

\paragraph{Description of $C$ and $\Gamma$.} There are explicit expressions for the deterministic coefficients $C$ and $\Gamma$ of theorems~\ref{thm:convergenceDonsker} and \ref{thm:convergenceembeddings} that are easy to derive but that we choose to skip here. Section~\ref{sec:studyCGamma} is dedicated to such computations giving in particular explicit expressions relating iterated occupation times with $C$ and $\Gamma$.

\paragraph{Generalized ergodic theorem for iterated occupation times.} Following the same type of proof as in section~\ref{sec:proofthm} based on excursion theory of Markov chain, one can prove the following asymptotic property.
\begin{proposition}
\label{prop:ergodic:occtime}
Let $(R_n)_{n\in\setN}$ be a Markov chain on a countable space $E$ such that it is irreducible and positive recurrent, with invariant probability $\nu$. Then, for any sequence $(u_1,\ldots,u_k)\in E_k$, 
\begin{equation}\label{eq:thm:ergodicocctime}
\frac{L_{u_1\ldots u_k}( (R_n)_{0\leq n<N} )}{N^k}
\xrightarrow[N\to\infty]{a.s., L^1} \frac{\nu(u_1)\nu(u_2)\ldots \nu(u_k)}{k!}
\end{equation}
\end{proposition}
This result looks uninteresting at first sight and corresponds to the discrete equivalent of the convergence of $\delta_{N^{-1}}\circ S\circ \iota_N(X)$ to a zero limit due to the centring of the increments and the law of large number. A more interesting asymptotic result for iterated occupation times consists in considering fluctuations around the a.s.~limit~\eqref{eq:thm:ergodicocctime}. This is done partially in section~\ref{sec:occtime:asymptotics}, even if it would be interesting to have a more general study of iterated occupation times.

\subsubsection{The general question of embeddings.}
Theorem~\ref{thm:convergenceDonsker} already encompasses a wide variety of discrete models but we formulate below a generalization of it by observing the two following facts. 

First, the definition~\ref{def:hiddenmarkovwalk} of hidden Markov walk requires only that $V$ is a semi-group in order to build $X_n=F_0\cdot F_1 \cdot\ldots\cdot F_{n-1}$ out of its increments $F_n$ and that it embeds nicely in $G^2(V)$ in order to formulate the theorem.

Besides, the choice of Donsker embedding is particular since it is one of the simplest embeddings one may consider on $V$ but one may choose more useful embeddings if the discrete model already has its own geometrical embedding (see for example the round-about model described in section~\ref{subsec:examples}).

We thus generalize the definition~\ref{def:hiddenmarkovwalk} to the case of general embeddings.

\begin{definition}
Let $E$ be a countable set and $V$ a finite-dimensional vector space. A hidden Markov path $((R_n,\mathbf{F}_n)_{n\in\setN},(X_t)_{t\in\setR_+})$ with regularity $\alpha$ is a process such that:
\begin{enumerate}[(i)]
\item the process $(R_n)_{n\in\setN}$ is a Markov chain on $E$;
\item the r.v. $\mathbf{F}_n$ have values in $C^{0,\alpha\text{-Höl}}([0,1],G_2(V))$
\item conditionally on the process $(R_n)$, the increments $\mathbf{F}_n$ are \emph{independent} and have marginal laws such that, for any Borel set $A$, any $n\in\setN$ and any $r\in E$, 
\begin{align}
\probc{ \mathbf{F}_n\in A}{\sigma(R)} &= \probc{ \mathbf{F}_n\in A}{R_n} 
\\
\probc{ \mathbf{F}_n\in A}{R_n=r} &= \probc{ \mathbf{F}_1\in A}{R_1=r}
\end{align}
and moreover $\mathbf{F}_n(0)=0_{G^2(V)}$ a.s.;
\item the process $(X_t)$ is obtained by concatenating the increments $\mathbf{F}_n$, i.e., for any $n\in\setN$, for any $t\in [n,n+1[$,
\begin{equation}
X_t =  \mathbf{F}_0(1)\cdot\ldots \mathbf{F}_{n-1}(1) \ldots \cdot \mathbf{F}_{n}(t-n)
\end{equation}
where $\cdot$ is the product in $G_2(V)$.
\end{enumerate}
\end{definition}
Given this generalized definition, theorem~\ref{thm:convergenceDonsker} admits the following generalization.
\begin{theorem}\label{thm:convergenceembeddings}
Let $E$ be a countable set and $V$ a finite-dimensional vector space. Let $r_0\in E$. Let $((R_n,\mathbf{F}_n)_{n\in\setN},(X_t)_{t\in\setR_+})$ be a hidden Markov path such that:
\begin{enumerate}[(i)]
\item $X_0=0_{G^2(V)}$ a.s. and $R_0=r_0$ a.s.
\item the Markov chain $(R_n)$ is irreducible and positive recurrent, with invariant probability $\nu$, and the first return time $T_1=\inf\{n\geq 1; R_n=R_0 \}$ has finite moments for all integer $p$:
\begin{equation}\label{eq:finitemomentsT1}
\Espi{r_0}{T_1^p} <\infty
\end{equation}
\item the increments $\mathbf{F}_n$ take values in $\bigcap_{\frac{1}{3}<\beta < \frac{1}{2}}C^{0,\beta\text{-Höl}}([0,1],G^2(V))$ and satisfy, for all $p\in\setN$, the bound
\begin{equation}
\label{eq:bounddistancealpha}
\sup_{r\in E} \sup_{\frac{1}{3}<\beta<\frac{1}{2}} \Espc{ d_\beta(\mathbf{F}_1,0_{G^2(V)})^p }{ R_1=r} <\infty
\end{equation}
\item the walk $(X_t)$ is centred, i.e. satisfies 
\begin{equation}
\label{eq:centering}
\sum_{r\in E} \nu(r) \Espc{ \pi_1(\mathbf{F}_1(1)) }{R_1=r} = 0_V
\end{equation}
where $\pi_1(u)$ is the component in $V$ of $u\in G^2(V)\subset V\oplus (V\otimes V)$.
\end{enumerate}
Then, for any $1/3<\alpha<1/2$, the sequence of processes $(\delta_{N^{-1/2}}(X_{Nt}))_{t\in[0,1]}$ converges in law in the rough path topology of $C^{0,\alpha\text{-Höl}}([0,1],G^2(V))$ to the enhanced Brownian motion (as defined in~\eqref{eq:enhancedanomalousBM}) with covariance matrix $C$ and area anomaly $\Gamma\in V\wedge V$ given in \eqref{eq:formulacovariancematrixC} and \eqref{eq:formulaareadrift}, for any $1/3<\alpha<1/2$. Moreover, the limit law does not depend on the choice of $r_0$.
\end{theorem}

In particular, one verifies easily that theorem~\ref{thm:convergenceDonsker} is a consequence of theorem~\ref{thm:convergenceembeddings} by choosing the linear interpolation: \begin{equation}\mathbf{F}_n(t) = \exp(tF_n) =  \left( t F_n, 0\right)\in V\oplus (V\wedge V) \simeq G^2(V)
\end{equation}
(see section~\ref{sec:introroughpath} for the exact definition of the Lie group $G^2(V)$). More generally, hypothesis~\eqref{eq:bounddistancealpha} is satisfied as soon as the embeddings are smooth or Lipschitz and the increments have finite moments.

As for theorem~\ref{thm:convergenceDonsker}, the area anomaly $\Gamma$ obtained in \eqref{eq:formulaforGammaocctimes} is related to the area of the walk $(X_n)$ covered during one excursion of the process $(R_n)$. There may be now two contributions: one of them is related to the fact that the $\ca{F}_n$ may contain an area drift (this is the new part due to the nonlinear embeddings) and one of them is related to the area produced by correlations during one excursion.

\subsection{Examples and heuristics}\label{subsec:examples}

We present in this section various models to illustrate the previous theorems. The explicit formulae for $C$ and $\Gamma$ are given below in equations~\eqref{eq:formulaforCocctimes} and \eqref{eq:formulaforGammaocctimes}.

\paragraph{A random walk in $\setC$ with rotating increments.} Let $\omega=e^{2i \pi/L}$ be a root of unity with $L\geq 3$. Let $(U_k)_{k\in\setN}$ be a sequence of i.i.d. real random variables with finite moments of all order. Identifying canonically $\setC$ with $\setR^2$, the process $(X_n)_{n\in\setN}$ defined by $X_0=0$ and, for any $n\geq 1$
\begin{equation*}
X_n= \sum_{k=0}^{n-1} \omega^k U_k
\end{equation*}
is a hidden Markov walk. Indeed, one may choose $E=\setZ/L\setZ$ and the deterministic dynamics $R_{n}= n \mod L$. The increment $F_n= \omega^{R_n} U_n$ depends only on $R_n$ and of r.v. independent of the process $(R_n)$. The first return time is constant $T_1=L$ and thus $\Espi{0}{T_1^p}=L^p <\infty$. Moreover one has the following computations.
\begin{align*}
\Espi{0}{X_{T_1}}&= \left(\sum_{k=0}^{L-1} \omega^k\right) \Esp{X_1} = 0
\\
C_{ij} &= \frac{1}{L}\Espi{0}{X^{(i)}_{T_1} X^{(j)}_{T_1}} = \frac{\Var(X_1)}{2}  \delta_{ij}
\\
\Gamma  &= \frac{\cos(\pi/L)}{4\sin(\pi/L)} \Esp{X_1}^2\begin{pmatrix} 0 & -1 \\ 1 & 0 \end{pmatrix}
\end{align*}
In the present case, the return time has an almost sure value and thus we may improve the proof of theorem~\ref{thm:convergenceembeddings} in order to relax the finite moment hypothesis on the $U_k$.

\paragraph{Spending time turning around.} A case described by theorem~\ref{thm:convergenceembeddings} but not by theorem~\ref{thm:convergenceDonsker} is given by the following construction. We fix $E=\{1,0\}$. If $R_n=1$, then a centred random vector $U_n$ is chosen and the path increment is the straight line:
\begin{equation}
\mathbf{F}_n(t) = (tU_n,0)
\end{equation}
If $R_n=0$, then the path increment is a circle $c(t)=r(\cos(2\pi t)-1,\sin(2\pi t))$. The covered area at time $1$ is thus $\pi r^2$ and $\mathbf{F}_n(1)=(0,\pi r^2)$ where $\setR^2\wedge\setR^2$ is identified to $\setR$. 

The process $R_n$ is a Markov chain with transition matrix \[ Q= \begin{pmatrix} 1-a & a \\ b & 1-b \end{pmatrix} \] in the basis $(1,0)$. Theorem~\ref{thm:convergenceembeddings} can be applied. An excursion corresponds to the sequence of states $1$ or  $100\ldots0$ where the number of $0$ is a geometric law (starting at $1$) and one obtains:
\begin{align*}
C &= \frac{1}{\Espi{1}{T_1}} \Espi{1}{X_{T_1}\otimes X_{T_1}} = \frac{b}{a+b} \Esp{U_1\otimes U_1}
\\
\Gamma &= \pi r^2  \frac{a}{a+b} 
\end{align*}
If $b$ goes to zero, one recovers the area bubbles mentioned in the introduction. If $a$ goes to $0$, there are no circles and one recovers the classical random walk.

\paragraph{Diamond and round-about models: the question of correlations and embeddings.} We introduce two other models that illustrate theorem~\ref{thm:convergenceembeddings} in figure~\ref{fig:diamondroundabout}. Both models have the same space $E=\{1,2,\ldots,8\}$. A value $r\in E$ corresponds to a unique type of edges on a $\setZ^2$-periodic graph. In both models, the edges have the same increments in the plane $\setR^2$ but differ by their embeddings: in the diamond model, all the embeddings are straight lines whereas, in the round-about model, part of the edges are circle arcs, which cover a non-zero area. 

At the end of an arrow, there is exactly two out-coming edges, one plain and one dashed. The dashed out-coming edge is chosen with probability $p$ and the plain one with probability $1-p$.

 If $p\geq 1/2$, the walk tends to be trapped in the diamonds/round-abouts and thus accumulates a covered area, which in the continuous limit, contributes to the area anomaly $\Gamma$. In the diamond model, all the contributions to $\Gamma$ are of these type. In the round-about model, there is an additional contribution to $\Gamma$ corresponding to the area covered by the circle arcs.

\begin{figure}
\begin{center}
\begin{tikzpicture}
\begin{scope}
\draw[->,dashed,ultra thick] (1,0) -- node [midway,above right] {$1$} (0,1) ;
\draw[->,dashed,ultra thick] (0,1) -- node [midway,above left] {$2$}(-1,0);
\draw[->,dashed,ultra thick] (-1,0) -- node [midway,below left] {$3$}(0,-1);
\draw[->,dashed,ultra thick] (0,-1) -- node [midway,below right] {$4$}(1,0);

\draw[->,ultra thick] (1,0.05) -- node [midway,above] {$5$}(2,0.05);
\draw[->,ultra thick] (-0.05,1) -- node [midway,right]{$7$}(-0.05,2);
\draw[->,ultra thick] (-1,-0.05) --  (-2,-0.05);
\draw[->,ultra thick] (0.05,-1) -- (0.05,-2);

\draw[->,ultra thick] (2,-0.05) --  node [midway,below] {$6$}(1,-0.05);
\draw[->,ultra thick] (0.05,2) -- node [midway,left] {$8$}(0.05,1);
\draw[->,ultra thick] (-2,0.05) -- (-1,0.05);
\draw[->,ultra thick] (-0.05,-2) -- (-0.05,-1);

\draw[dashed,ultra thick] (2,0) -- (2.5,-0.5);
\draw[->,dashed,ultra thick] (2.5,0.5) -- (2,0);

\draw[dashed,ultra thick] (0,2) -- (0.5,2.5);
\draw[->,dashed,ultra thick] (-0.5,2.5) -- (0,2);

\draw[dashed,ultra thick] (-2,0) -- (-2.5,0.5);
\draw[->,dashed,ultra thick] (-2.5,-0.5) -- (-2,0);

\draw[dashed,ultra thick] (0,-2) -- (-0.5,-2.5);
\draw[->,dashed,ultra thick] (0.5,-2.5) -- (0,-2);

\draw[dotted, thick] (-2.5,-2.5)--(2.5,-2.5)-- (2.5,2.5) -- (-2.5,2.5) -- (-2.5,-2.5);
\end{scope}

\begin{scope}[xshift=6cm]
\draw[->,dashed,ultra thick] (1,0) arc (0:85:1);
\draw[->,dashed,ultra thick] (0,1) arc (90:180:1);
\draw[->,dashed,ultra thick] (-1,0) arc (180:270:1);
\draw[->,dashed,ultra thick] (0,-1) arc (270:360:1);
\node at (0.9,0.9) {$1$};
\node at (-0.9,0.9) {$2$};
\node at (-0.9,-0.9) {$3$};
\node at (0.9,-0.9) {$4$};
\node at (1.5,0.3) {$5$};
\node at (1.5,-0.3) {$6$};
\node at (0.3,1.5) {$7$};
\node at (-0.3,1.5) {$8$};

\draw[->,ultra thick] (1,0.05) -- (2,0.05);
\draw[->,ultra thick] (-0.05,1) -- (-0.05,2);
\draw[->,ultra thick] (-1,-0.05) -- (-2,-0.05);
\draw[->,ultra thick] (0.05,-1) -- (0.05,-2);

\draw[->,ultra thick] (2,-0.05) -- (1,-0.05);
\draw[->,ultra thick] (0.05,2) -- (0.05,1);
\draw[->,ultra thick] (-2,0.05) -- (-1,0.05);
\draw[->,ultra thick] (-0.05,-2) -- (-0.05,-1);

\draw[dashed,ultra thick] (2,0) arc (180:240:1);
\draw[->,dashed,ultra thick] (2.5,0.86) arc (120:180:1);

\draw[dashed,ultra thick] (0,2) arc (270:330:1);
\draw[->>,dashed,ultra thick] (-0.86,2.5) arc (210:270:1);

\draw[dashed,ultra thick] (-2,0) arc (0:60:1);
\draw[->,dashed,ultra thick] (-2.5,-0.86) arc (-60:0:1);

\draw[dashed,ultra thick] (0,-2) arc (90:150:1);
\draw[->,dashed,ultra thick] (0.86,-2.5) arc (30:90:1);

\draw[dotted, thick] (-2.5,-2.5)--(2.5,-2.5)-- (2.5,2.5) -- (-2.5,2.5) -- (-2.5,-2.5);
\end{scope}
\end{tikzpicture}
\end{center}
\caption{\label{fig:diamondroundabout}
Diamond (left) and round-about (right) models on $\setZ^2$-periodic patterns. At each vertex, there is exactly two incoming and two out-coming edges, one plain and one dashed. The dashed out-coming edge is chosen with probability $p$ and the plain one with probability $1-p$. The space state $E$ corresponds to the types of arrow. The difference between the two cases corresponds to the curved dashed arrow: in the roundabout model, the increment is the same but with an additional covered area.}
\end{figure}

\section{Mathematical tools and proofs.}
\subsection{The rough path space and its topology}
\label{sec:introroughpath}

All our notations and definitions follow closely the ones introduced in \citep{Breuillardetal} and \citep{FrizVicBook} and thus we sketch only the tools needed for the proofs. The space $G^2(V)$ is the subset of $T_1^{(2)}(V) = V\oplus (V\otimes V)$ of elements $(v,M)$ such that there exists a smooth path $x: [0,1]\to V$ whose signature $S(x)(1)$ is equal to $(v,M)$. It is easy to see that $G^2(V)$ is a subgroup isomorphic to $V\oplus V\wedge V$ since the symmetric part of $M$ has to be equal to $(1/2)v\otimes v$ and thus can be skipped from the description.
The group law is defined as
\[
(a,A)\cdot(b,B) = (a+b,A+B+(1/2)(a\otimes b-b\otimes a))
\]
and the inverse is given by $(a,A)^{-1}=(-a,-A)$. We define the two canonical projections:
\begin{align*}
\pi_1 : G^2(V) & \to V
& 
\pi_2 : G^2(V) & \to V \wedge V
\\
(a,A) & \mapsto a
&
(a,A) & \mapsto A
\end{align*}

The Carnot-Caratheodory norm $\norm{u}$ of an element $u$ is the infimum of the lengths\footnote{$V$ is assumed to be Euclidean.} of smooth paths $x$ such that $S(x)(1)=u$ and it induces a distance on $G^2(V)$ through $d(u_1,u_2)=\norm{u_1\cdot u_2^{-1}}$, making $G^2(V)$ a geodesic space.

Given two smooth paths $x_1,x_2:[0,1]\to V$ such that $x_1(0)=x_2(0)=0_V$, we introduce, for any $\alpha\in(1/3,1/2)$, the distance:
\begin{equation}
\label{eq:distancealpha}
d_{\alpha}(x_1,x_2) = 
\sup_{(s,t)\in [0,1]^2}
\frac{ \norm{S_2(x_1)(s,t)\cdot S_2(x_2)(s,t)^{-1}}}{|t-s|^\alpha}
\end{equation}
where $S_2(x)(s,t)=S_2(x)(t)\cdot S_2(x)(s)^{-1}$ is the $G^2(V)$-valued increment of the signature.

The $d_\alpha$-closure of the set of signatures $(S(x))$ of smooth paths $x$ at finite distance from the signature of the constant zero path is the Polish space (for $d_\alpha$)
\begin{equation}
C^{0,\alpha\text{-Höl}}([0,1],G^2(V))
\end{equation}

In practice, the study of convergence in law of processes in the topology of this rough path space is made easier by two useful tools: equivalence of norms and Kolmogorov-Centsov tightness criterion. First, there is equivalence of norms on $G^2(V)$ between the Carnot-Caratheodory norm $\norm{(a,A)}$ defined through geodesics and the norm $\norm{(a,A)}'$ defined by:
\begin{equation}
\label{eq:defnormprime}
\norm{(a,A)}'= \sum_{i=1}^d |\pi_1^{(i)}(a)| + \sum_{1\leq i<j\leq d} |\pi_2^{(ij)}(A)|^{1/2}
\end{equation}
where $\pi_1^{(i)}(a)$ and $\pi_2^{(ij)}(A)$ are the components of $a$ and $A$ in a given basis. Bounding $d_\alpha(x_1,x_2)$ thus only requires suitable bounds on every coefficient of the components of the paths.

All the tightness criteria required by our theorem deal with $1/3< \alpha < 1/2$ and thus, following \citep{Breuillardetal}, tightness of a sequence of processes $((X^{(N)}(t))_{t\in[0,1]})_{N\geq 1}$ in a fixed $C^{0,\alpha\text{-Höl}}([0,1],G^2(V))$ requires only that there exists $p\geq 1$ such that $\alpha\leq (2p-1)/(4p)$ and there exists $C$ such that for any $s,t\in[0,1]$, it holds:
\begin{equation}\label{eq:kolmorogovtightness}
\sup_{N\geq 1} \Esp{ \norm{ X^{(N)}(t)\cdot X^{(N)}(s)^{-1} }^{4p} } \leq C |t-s|^{2p}
\end{equation}

\subsection{Proof of the convergence theorem and formulae for the anomalous area drift}
\label{sec:proofthm}

In this section, we fix a hidden Markov path  $((R_n,\mathbf{F}_n)_{n\in\setN},(X_t)_{t\in\setR_+})$ in $G^2(V)$ satisfying the hypothesis of theorem~\ref{thm:convergenceembeddings}. The proof of theorem~\ref{thm:convergenceembeddings} relies on the following steps: 
\begin{itemize}
\item we cut the trajectories of $(R_n)$ into excursions;
\item we then study the convergence of the accelerated geodesic interpolation $(\ha{X}_t)$ of the process $(X_t)$ between two successive return times of $(R_n)$;
\item we compare the finite-dimensional marginals of the two processes $\ha{X}$ and $X$ in the limit;
\item we prove the tightness of the sequence of processes $(X^{(N)})_N$.
\end{itemize}

\paragraph{Cutting into independent excursions.} The proof of theorem~\ref{thm:convergenceembeddings} relies on the division of the process $((R_n,\mathbf{F}_n)_{n\in\setN},(X_t)_{t\in\setR_+})$ into time windows $[T_k,T_{k+1})$ corresponding to excursions of the Markov process $(R_n)$ (see \citep{bookNorris} for a general theory of excursions of Markov processes).

\begin{proposition}
\label{prop:excursionincrements}
Let $(T_k)_{k\in\setN}$ be the sequence of excursion times defined by:
\begin{equation}
\begin{cases}
T_0 &= 0
\\
T_{k+1} &= \inf\left\{ n>T_k ; R_{n}=R_{T_k} \right\}, 
\qquad k\in\setN
\end{cases}
\end{equation}
Then, for any $r_0\in E$, conditionally on $\{R_0=r_0\}$, the r.v. $(\ha{F}_k)_{k\in\setN_1}$ defined by
\begin{equation}\label{eq:defexcursionincrements}
\ha{F}_k = \mathbf{F}_{T_{k-1}}(1)\cdot\mathbf{F}_{T_{k-1}+1}(1)\cdot \ldots \cdot\mathbf{F}_{T_{k}-1}(1)
\end{equation}
is a sequence of independent and identically distributed $G^2(V)$-valued random variables.
\end{proposition}
\begin{proof}
The recurrence of $(R_n)$ implies that all the $T_k$ are finite a.s. and, from the general theory of discrete Markov processes, the excursions of the process $(R_n)$ are independent and identically distributed. Each $\ha{F}_k$ is a product of r.v. indexed by times belonging to the same excursion of the process $(R_n)$ and the hidden Markov structure implies the result.
\end{proof}

\begin{property}[moments of $\ha{F}_1$]
\label{prop:momentsofhatF}
Under the hypotheses of theorem~\eqref{thm:convergenceembeddings}, the r.v.~$(\ha{F}_k)_{k\in\setN}$ are independent and satisfy:
\begin{enumerate}[(i)]
\item for all $p\in\setN$ and all $r_0\in E$, $\Espi{r_0}{\norm{\ha{F}_1}^p} < \infty$
\item the projection on $V$ are centred and have a finite covariance
\begin{align}\label{eq:formulacovariancematrixC}
\Espi{r_0}{\pi_1(\ha{F}_1)} &= 0_V 
&
\Espi{r_0}{\pi_1^{(i)}(\ha{F}_1)\pi_1^{(j)}(\ha{F}_1) } &= C_{ij} \Espi{r_0}{T_1} < \infty
\end{align}
where $C_{ij}$ is symmetric and does not depend on the choice of $r_0$.
\item the expectation of the second level is given by
\begin{equation}\label{eq:formulaareadrift}
\Espi{r_0}{ \pi_2^{(ij)}(\ha{F}_1) } = \Gamma_{ij} \Espi{r_0}{T_1}
\end{equation}
where $\Gamma_{ij}$ is antisymmetric and does not depend on $r_0$.
\end{enumerate}
\end{property}
\begin{proof}
Let $p\in\setN$ be a fixed integer. The first point uses sub-additivity of the norm and the bound $(a_1+\ldots+a_n)^p \leq n^p (a_1^p+\ldots+a_n^p)$ for positive numbers:
\begin{align*}
\norm{\ha{F}_1}^p
& \leq T_1^p \sum_{0\leq k< T_1} \norm{\mathbf{F}_k(1)}^p
\\
\Espci{r_0}{ \norm{\ha{F}_1}^p}{\sigma(R)} 
& \leq T_1^p \sum_{0\leq k< T_1} \Espci{r_0}{\norm{\mathbf{F}_k(1)}^p}{\sigma(R)} \leq T_1^{p+1} C
\end{align*}
using hypotheses~\eqref{eq:bounddistancealpha} and \eqref{eq:finitemomentsT1}. 

The second point uses the classical representation property of invariant measure as marginal of the excursion measure of additive functionals: for any $\nu$-integrable function $f:E\to\setR$, it holds 
\begin{equation}
\sum_{r\in E} f(r) \nu(r)  = \frac{1}{\Espi{r_0}{T_1}} \Espi{r_0}{ \sum_{0\leq k< T_1} f(R_k)}
\end{equation}
for any $r_0\in E$. In the present case, we apply this formula to \begin{align*}
\Espi{r_0}{\pi_1(\ha{F}_1)} = \Espi{r_0}{\sum_{0\leq k<T_1} \Espci{r_0}{\pi_1(\mathbf{F}_k(1))}{R_k}}
\end{align*}
and use the centring hypothesis~\eqref{eq:centering}. $C$ and $\Gamma$ do not involve additive functionals (see below section~\ref{sec:studyCGamma} for more information and explicit formulae), however the independence with respect with $r_0$ can be proved in the same way as the previous property.
\end{proof}

\paragraph{Convergence of the excursion-geodesic extracted process.}

Out of the independent $G^2(V)$-valued r.v.~$(\ha{F}_k)$, we follow \citep{Breuillardetal} and build the geodesic-interpolated processes $(\ha{X}_{t}))_{t\in \setR_+}$ defined by
\begin{equation}
\ha{X}_{t} = \ha{F}_1\cdot\ldots\cdot \ha{F}_{\IntP{t}} \cdot g(\ha{F}_{\IntP{t}+1},\FracP{t})
\end{equation}
where $g: G^2(V) \times [0,1]\to G^2(V)$ is defined such that $g(u,\cdot)$ is the geodesic in $G^2(V)$ joining $0_{G^2(V)}$ and $u$. For any real number $t$, $\IntP{t}$ and $\FracP{t}$ are respectively its integer and fractional parts.

The relation with the initial process $X_{t}$ is that,
for any $k\in\setN$,
\begin{equation}\label{eq:acceleratedhatprocess}
\ha{X}_{k} = X_{T_k}
\end{equation}
However, the sequence of processes $(\delta_{N^{-1/2}}(\ha{X}_{Nt}))$ is such that a corrected version of the results of \citep{Breuillardetal} can be applied and provide the following theorem and the following lemma. By the word \emph{corrected}, we mean that theorem~3 of \citep{Breuillardetal} should either require also the centering of the $\pi_2(\mathbf{\xi}_i)$ in order to have no drift area or, if not, include an area drift 
$\Gamma$ in the enhanced Brownian motion at the limit. One may be convinced for example by considering increments $(\xi_i,a)$ where $\xi_i$ is an random increment in $V$ and $a$ is a constant in $V\wedge V$. Then $\mathbf{W}_{k/n}^{(n)}$ is given by $\delta_{n^{-1/2}}\left( e^{\xi_1}\otimes\ldots\otimes e^{\xi_k}\right)\otimes (0,ka/n)$ since the last term is central in $G^2(V)$ and converges a.s. to the limit process $(0,ta)$ (using Slutsky's lemma, one has convergence to the anomalous enhanced Brownian motion). The proof of theorem~3 of \citep{Breuillardetal} remains unchanged: the use of Stroock-Varadhan theorem is still valid but identifies a non-zero additional drift term in $V\otimes V$. The tightness criterion remains the same up to recentering of the area, which is costless. An alternative way of identifying this drift term is present in \citep{ChevyrevLevy} in a more general context and with a precise description of all the terms.

\begin{theorem}[from \citep{Breuillardetal}]\label{thm:excursionextraction}
Let $((R_n,\mathbf{F}_n)_{n\in\setN},(X_t)_{t\in\setR_+})$ be a hidden Markov path in $G^2(V)$ satisfying the hypothesis of theorem~\ref{thm:convergenceembeddings} and let $1/3<\alpha<1/2$. Let $\beta=1/\Espi{r_0}{T_1}$. The sequence of processes $((\delta_{N^{-1/2}}(\ha{X}_{N\beta t})_{t\in[0,1]})_{N\in\setN^*}$ converges in distribution in the space $C^{0,\alpha\text{-Höl}}([0,1],G^2(V))$ to the enhanced Brownian motion with covariance matrix $C$ and area anomaly $\Gamma$ given in property~\ref{prop:momentsofhatF}. Moreover, the limit law does not depend on $r_0$.
\end{theorem}

We emphasize that the dependence on $r_0$ in the construction of $\ha{X}_k$ is due to the construction of the excursions of the process $(R_n)$ and of their length $T_1$: the independence of $C$ and $\Gamma$ in property~\ref{prop:momentsofhatF} requires a normalization by $\Espi{r_0}{T_1}$, which is included in the previous theorem as a slow-down of the time scale. This can be also seen in eq.~\eqref{eq:acceleratedhatprocess}, using the a.s. asymptotic equivalent $T_k\simeq \Esp{T_1}k$.

\paragraph{Useful bounds on products of increments $\ha{F}_n$.}

\begin{lemma}[from \citep{Breuillardetal}]
For all integer $p\geq 1$, there exists $C'_p>0$ such that, for any integers $n<m$, the following bound holds:
\begin{equation}
\Esp{
\norm{ \ha{F}_n\cdot \ha{F}_{n+1}\cdot\ldots\cdot\ha{F}_{m-1} }^{4p}
} \leq C_p (m-n)^{2p} 
\end{equation}
\end{lemma}
This lemma inherited directly from~\citep{Breuillardetal} corresponds to equation~$(*)$ of \citep{Breuillardetal} and its proof relies on the centering of the r.v. $\pi_1^{(i)}(\ha{F}_n)$ and a equivalence of $\norm{\cdot}$ with a norm which is a polynomial in the components of $\ha{F}_n$.

However, this lemma is not precise enough for our purposes because of the fluctuations of the time scale between $\ha{X}_{N\beta t}$ and the process $X_{Nt}$ (more precisely $T_k/k \to \Esp{T_1}$ only asymptotically) and we need the following improved \emph{maximal} version of the previous lemma.

\begin{lemma}
\label{lemma:boundexcursionfortightness}
For all integer $p\geq 1$, there exist $C'_p>0$  such that, for any integers $n<m$, the following bounds hold:
\begin{align}
\label{eq:maxineq1}
\Esp{
\sup_{n\leq k \leq m-1} \norm{ \ha{F}_n\cdot\ha{F}_{n+1}\cdot\ldots\cdot\ha{F}_{k} }^{4p}
} & \leq C'_p (m-n)^{2p} 
\\
\label{eq:maxineq2}
\Esp{
\sup_{n\leq k\leq l \leq m-1} \norm{ \ha{F}_k\cdot\ha{F}_{k+1}\cdot\ldots\cdot\ha{F}_{l} }^{4p}
} & \leq 2^{4p}C'_p (m-n)^{2p} 
\end{align}
\end{lemma}
\begin{proof}
We start with the proof of \eqref{eq:maxineq1}. Using the distance \eqref{eq:defnormprime}, it is enough to prove the existence of constants $A_{i}$ and $B_{ij}$ such that:
\begin{subequations}
\begin{align}
\Esp{\sup_{n\leq k \leq m-1} \left|\pi_1^{(i)}( \ha{F}_n\cdot\ha{F}_{n+1}\cdot\ldots\cdot\ha{F}_{k} )\right|^{4p}} 
& \leq A_i (m-n)^{2p}
\label{eq:pi1bound}
\\
\Esp{\sup_{n\leq k \leq m-1} \left|\pi_2^{(ij)}( \ha{F}_n\cdot\ha{F}_{n+1}\cdot\ldots\cdot\ha{F}_{k} )\right|^{2p}} 
& \leq B_{ij} (m-n)^{2p}
\label{eq:pi2bound}
\end{align}
\end{subequations}
for any $1\leq i,j \leq d$. By definition of $G^2(V)$, we have:
\begin{align*}
\pi_1^{(i)}( \ha{F}_n\cdot\ha{F}_{n+1}\cdot\ldots\cdot\ha{F}_{k} )
= &\sum_{l=n}^{k} \pi_1^{(i)}(\ha{F}_l)
\\
\pi_2^{(ij)}( \ha{F}_n\cdot\ha{F}_{n+1}\cdot\ldots\cdot\ha{F}_{k} )
= & \sum_{l=n}^{k} \pi_2^{(ij)}(\ha{F}_l) 
\\
& + \sum_{n\leq l_1 < l_2 \leq k} \left( \pi_1^{(i)}(\ha{F}_{l_1})\pi_1^{(j)}(\ha{F}_{l_2}) - 
\pi_1^{(j)}(\ha{F}_{l_1})\pi_1^{(i)}(\ha{F}_{l_2}) 
\right)
\end{align*}
The r.v.~$\pi_1^{(i)}(\ha{F}_l)$ are i.i.d.~centred random variable and thus the sequence $(M^{(1)}_k)_{k\geq n}$ of r.v. $M^{(1)}_k=\sum_{l=n}^k \pi_1^{(i)}(\ha{F}_l)$ is a martingale  and Doob's maximal inequality gives the first bound \eqref{eq:pi1bound} since $\Esp{(M^{(1)}_{m-1})^{4p}} \leq A'_i (m-n)^{2p}$. 

We introduce now $a^{(ij)}=\Esp{\pi_2^{(ij)}(\ha{F}_1)}$ and observe that, if $M^{(2)}_k=\pi_{2}^{(ij)}(\ha{F}_n\cdot\ldots\cdot\ha{F}_k)-(k-n)a^{(ij)}$, $(M^{(2)}_k)_{k\geq n}$ is also a martingale since both $\left(\sum_{l=n}^k (\pi_2^{(ij)}(\ha{F}_l)-a^{(ij)})\right)_{k\geq n}$ and \[
\left(\sum_{n\leq l_1 < l_2 \leq k} \left( \pi_1^{(i)}(\ha{F}_{l_1})\pi_1^{(j)}(\ha{F}_{l_2}) - 
\pi_1^{(j)}(\ha{F}_{l_1})\pi_1^{(i)}(\ha{F}_{l_2}) 
\right) \right)_{k\geq n}
\]
are martingales. Moreover it holds:
\[\left|\pi_2^{(ij)}( \ha{F}_\cdot\ldots\cdot\ha{F}_{k} )\right|^{2p} \leq 2^{2p} \left( |a^{(ij)}|^{2p} (k-n)^{2p} + M_k^{2p} \right) \leq C_1 (m-n)^{2p} + M_k^{2p}
\]
Using again Doob's maximal inequality, we are left to prove
\begin{equation}
\Esp{(M^{(2)}_{m-1})^{2p}} \leq C_2(m-n)^{2p}
\end{equation}
which can be done by direct computation. Hence we obtain eq.~\eqref{eq:pi2bound}.

The proof of \eqref{eq:maxineq2} is obtained by the norm subadditivity through the following inequality for $n\leq k \leq l \leq m$
\begin{align*}
\norm{\ha{F}_k\cdot\ldots\cdot \ha{F}_l}^{4p}&= \norm{(\ha{F}_n\cdot\ldots\cdot\ha{F}_{k-1})^{-1}\cdot(\ha{F}_n\cdot\ldots\cdot\ha{F}_{l})}^{4p}
\\
&\leq 2^{4p} \left(\norm{\ha{F}_n\cdot\ldots\cdot\ha{F}_{k-1}}^{4p} + \norm{\ha{F}_n\cdot\ldots\cdot\ha{F}_{l}}^{4p}\right)
\end{align*}
Using then \eqref{eq:maxineq1} gives directly \eqref{eq:maxineq2}.
\end{proof}

\paragraph{Useful bounds on the renewal process $(T_k)_{k\in\setN}$.}

The sequence of excursion times $(T_k)_{k\in\setN}$ is a renewal process on $\setN$ (see \citep{Asmussen} for a good introduction). We introduce, for any $x\in\setR_+$, the integer-valued random variable $K(x)$ defined by
\[
K(x) = \max\left\{ n\in\setN ; T_n \leq x \right\}
\]
This is the unique integer $K(x)$ such that $T_{K(x)}\leq x< T_{K(x)+1}$. The proofs below require deviation estimates of $K(x)-\beta x$ with $\beta=1/\Esp{T_1}$ as well as moment estimates on increments $T_{K(x)+1}-T_{K(x)}$. We emphasize on the \emph{inspection paradox}: although $T_{n+1}-T_n$ has the same law as $T_1$, $T_{K(x)+1}-T_{K(x)}$ is stochastically larger than $T_1$. However, in our case, the following property still holds:
\begin{proposition}\label{prop:inspectionmoments}
Let $(T_k)$ be a renewal process such that the $(T_{n+1}-T_{n})_{n\in\setN}$ are i.i.d. and there exists $P\in \setN$ such that, for all $p\leq P$, $\Esp{T_1^p}<\infty$. Then, for all $p\leq P-1$, there exists $C_p>0$ such that, for all $x\in\setR_+$, 
\[
\Esp{  \left( T_{K(x)+1}-T_{K(x)}\right)^p } \leq C_p <\infty
\]
\end{proposition} 
\begin{proof}
It is sufficient to prove the result for $x$ integer since $K(\IntP{x})=K(x)$. We introduce the delay $D_j=T_{j}-T_{j-1}$ for $j\geq 1$ and we now decompose the event $\{ D_{K(x)+1}=n \}$ along the values of $K(x)$:
\begin{align*}
\{ D_{K(x)+1}=n \} &= \bigcup_{j\in\setN} \{ K(x)=j\;\text{and}\; D_{j+1}=n \} 
\\
&= \bigsqcup_{j\in\setN} \{ T_j \leq x < T_j+n \;\text{and}\; D_{j+1}=n\}
\\
&= \bigsqcup_{x-n < k \leq x} \bigsqcup_{j\in\setN}
\{ T_j =k \;\text{and}\; D_{j+1}=n\}
\end{align*}
We now obtain 
\[
\prob{ D_{K(x)+1}=n } = \sum_{x-n< k \leq x}\sum_{j\in \setN} \prob{
 T_j =k \;\text{and}\; D_{j+1}=n}
\]
For fixed $j$, the r.v. $T_j$ and $D_{j+1}$ are independent by construction and $D_{j+1}$ has the same law as $T_1$. We thus obtain
\begin{align*}
\prob{ D_{K(x)+1}=n } &= \prob{T_1=n}\sum_{x-n< k \leq x}\sum_{j\in \setN} \prob{T_j =k}
\\
&\leq \prob{T_1=n} \sum_{x-n< k \leq x} 1 \leq n\prob{T_1=n}
\end{align*}
The expected uniform finite bounds on the moments of $T_{K(x)+1}-T_{K(x)}$ is then an easy consequence.
%
%
%
\end{proof}
We will also need the following deviation estimation on the r.v. $K(x)$ for the proof of lemma~\ref{lemma:fdconvergence}.
\begin{proposition}\label{prop:deviationK}
Let $(T_k)$ be a renewal process such that the $(T_{n+1}-T_{n})_{n\in\setN}$ are i.i.d. and, there exists $p\in\setN$, $\Esp{T_1^p}<\infty$. Let $\beta=1/\Esp{T_1}$. Then, for all real $a>1$ and $x>0$, it holds:
\begin{subequations}
\begin{align}
\prob{ K(x)-\beta x \geq a } &\leq C_p \frac{(\beta x+a)^{p/2}}{(a-1)^p} \label{eq:Kdeviation1}
\\
\prob{ K(x)-\beta x \leq -a } & \leq \begin{cases} 
C_p\frac{(\beta x-a +1)^{p/2}}{(a-1)^p}
& \text{for $a\leq\beta x$} \\
0 & \text{for $a>\beta x$}
\end{cases}
\label{eq:Kdeviation2}
\end{align}
\end{subequations}
\end{proposition}
\begin{proof}
We first relate the event $\{|K(x)-\beta x| \geq a\}$ to events related to the $(T_n)$ using the equality $\{K(x) \geq n\}=\{ T_{n} \leq x\}$:
\begin{align*}
\prob{K(x)-\beta x \geq a }
& =  \prob{K(x) \geq \IntP{\beta x + a} } 
\leq  \prob{ T_{\IntP{\beta x+a}} \leq x } 
\\
&\leq 
\prob{ T_{\IntP{\beta x+a}}-\Esp{T_1}\IntP{\beta x+a} \leq -(a-1)\Esp{T_1} } 
\end{align*}
As a sum of independent centered r.v. with finite moments, it holds 
$\Esp{|T_{n}-n\Esp{T_1}|^p} \leq A_p n^{p/2}$. Using now Markov inequality produces \eqref{eq:Kdeviation1}. The second inequality is obtained in the same way by considering complementary events:
\[
\prob{K(x)-\beta x \leq -a } \leq \prob{K(x)\leq \IntP{\beta x-a} }
=  \prob{ T_{\IntP{\beta x-a}+1} > x}
\]
\end{proof}

\paragraph{Comparison of the finite-dimensional marginals of the two processes.} 

\begin{lemma}\label{lemma:fdconvergence}
Let $\beta=1/\Esp{T_1}\in (0,+\infty)$. Let $t\in [0,1]$. For any $\epsilon>0$, it holds:
\begin{equation}
\prob{ \norm{ \delta_{N^{-1/2}}(\ha{X}_{N\beta t})^{-1}\cdot \delta_{N^{-1/2}}(X_{Nt})} > \epsilon } \xrightarrow{N\to\infty} 0
\end{equation}
\end{lemma}
Using the same argument repeatedly, the same convergence in probability to $0$ can be extended to a family $(t_1,\ldots,t_m)$ and thus, by Slutsky's lemma, it shows that both processes have the same finite-dimensional marginal laws.

\begin{proof}
The idea behind the lemma comes from eq.~\eqref{eq:acceleratedhatprocess}: the stopping times $T_k$ are asymptotically equivalent to $k \Esp{T_1}$. The differences between the two processes have two origins: first, the difference between $T_k$ and $k\Esp{T_1}$ and second, the difference between the geodesic interpolation and the stochastic path $X_{t}$.

The two r.v. $\ha{X}_{Nt\beta}$ and $X_{Nt}$ are given by:
\begin{align*}
\ha{X}_{N\beta t}
&= 
\ha{F}_1\cdot\ldots\cdot\ha{F}_{\IntP{N\beta t}}\cdot g(\ha{F}_{\IntP{Nt\beta}+1},\FracP{Nt\beta})
\\
X_{Nt}
&= \ha{F}_1\cdot\ldots\cdot \ha{F}_{K(\IntP{Nt})}\cdot \mathbf{F}_{T_{K(Nt)}}(1)\cdot\ldots\cdot \mathbf{F}_{\IntP{Nt}-1}(1)\cdot \mathbf{F}_{\IntP{Nt}}(\FracP{Nt})
\end{align*}
where $K(u)$ is the unique integer such that
\[
T_{K(u)} \leq u < T_{K(u)+1}.
\]
We then obtain the bound:
\begin{align*}
\norm{\ha{X}_{N\beta t}^{-1} \cdot X_{Nt}}^{q}
\leq\;& 
3^q \Big( \norm{g(\ha{F}_{\IntP{N\beta t}+1},\FracP{Nt\beta})}^q
\\
&+ \norm{ \ha{F}_{\min(K(Nt),\IntP{Nt\beta})+1}\cdot \ldots 
\cdot\ha{F}_{\max(K(Nt),\IntP{Nt\beta})}}^q
\\
&+ \norm{ \mathbf{F}_{T_{K(Nt)}}(1)\cdot\ldots \cdot\mathbf{F}_{\IntP{Nt}-1}(1) \cdot\mathbf{F}_{\IntP{Nt}}(\FracP{Nt})}^q
\Big)
\end{align*}
The first norm is trivially bounded by $\norm{\ha{F}_{\IntP{N\beta t}+1}}^q$, which has a finite expectation from property~\ref{prop:momentsofhatF}. The third norm is also easy to bound in expectation:
\begin{align*}
&\Espc{
\norm{ \mathbf{F}_{T_{K(Nt)}}(1)\cdot\ldots \cdot\mathbf{F}_{\IntP{Nt}-1}(1)\cdot \mathbf{F}_{\IntP{Nt}}(\FracP{Nt})}^q
}{\sigma(R)} 
\\
&\leq (\IntP{Nt}-T_{K(Nt)})^{q} \sum_{k=T_{K(Nt)}}^{\IntP{Nt}} \Espc{ d_\alpha(\mathbf{F}_{k},0_{G^2(V)})^q}{ R_k} \\
&\leq A (\IntP{Nt}-T_{K(Nt)})^{q+1}
\leq A (T_{K(Nt)+1}-T_{K(Nt)})^{q+1}
\end{align*}
and thus the expectation of the $q$-th power of the norm is finite from proposition~\ref{prop:inspectionmoments} and hypothesis~\ref{eq:finitemomentsT1}.

In order to bound from above the expectation of the second norm, we introduce the sequence $u_N=\IntP{N^{\gamma}}$ for some real $\gamma>0$ chosen below. We then obtain, using $q=4p$ and lemma~\ref{lemma:boundexcursionfortightness}:
\begin{align*}
&\Esp{ \norm{ \ha{F}_{\min(K(Nt),\IntP{Nt\beta})-1} \cdot\ldots\cdot 
\ha{F}_{\max(K(Nt),\IntP{Nt\beta})}}^{4p} \indic{|K(Nt)-\IntP{Nt\beta}|\leq u_N }} 
\\
&\leq \Esp{ \sup_{\IntP{N\beta t}-u_N\leq k \leq l \leq  \IntP{N\beta t}+u_N} \norm{ \ha{F}_{k} \cdot\ldots\cdot 
\ha{F}_{l}}^{4p} } \leq K_1 u_N^{2p}
\end{align*}
The event $|K(Nt)-\IntP{N\beta t}|> u_N$ has a small probability for large $N$. Moreover the previous norm contains at most $\IntP{Nt}$ terms since $K(Nt)\leq \IntP{Nt}$ and $\beta > 1$. We use lemma~\ref{lemma:boundexcursionfortightness} and Cauchy-Schwarz inequality to obtain
\begin{align*}
J(u_N)= & \Esp{ \norm{ \ha{F}_{\min(K(Nt),\IntP{Nt\beta})} \cdot\ldots \cdot
\ha{F}_{\max(K(Nt),\IntP{Nt\beta})-1}}^{4p} \indic{|K(Nt)-\IntP{Nt\beta}|> u_N }} 
\\
\leq & \Esp{ \sup_{1\leq k \leq l \leq \IntP{Nt}} \norm{ \ha{F}_{k} \cdot\ldots \cdot
\ha{F}_{l}}^{4p} \indic{|K(Nt)-\IntP{Nt\beta}|> u_N }} 
\\
\leq & K_2 N^{2p}\,\prob{ |K(Nt)-\IntP{Nt\beta}|>u_N}^{1/2}
\end{align*}
Proposition~\ref{prop:deviationK} for $p=2$ immediately gives the following bound for $N$ large enough:
\[
J(u_N)\leq K_4 \frac{N^{2p+1/2}}{u_N}
\]
Collecting all the results with the dilation $\delta_{N^{-1/2}}$ gives for $4p+3\leq r$
\begin{align*}
\Esp{ \norm{\delta_{N^{-1/2}}(\ha{X}_{N\beta t})^{-1}\cdot\delta_{N^{-1/2}}(X_{Nt})}^{4p}}
\leq \frac{A_1}{N^{2p}} + \frac{K_1 u_N^{2p}}{N^{2p}} +
K_4\frac{N^{1/2}}{u_N}
\end{align*}
Any choice $1/2<\gamma < 1$ implies that the expectation tends to zero, hence the convergence in probability.
\end{proof}

\paragraph{Tightness of the initial process.}

\begin{lemma}[tightness]
Under the hypotheses of theorem~\ref{thm:convergenceembeddings}, the tightness criterion \eqref{eq:kolmorogovtightness} holds for the sequence of processes $(\delta_{N^{-1/2}}(X_{Nt}))_{t\in[0,1]}$.
\end{lemma}
\begin{proof}
The proof is similar to the one of the previous lemma. We fix $s<t$. If $\IntP{Ns}=\IntP{Nt}=j$, then 
\[
\norm{X_{Ns}(s)^{-1}\cdot X_{Nt}(t)}^{4p} = \norm{\mathbf{F}_{j}(Nt-Ns)}^{4p}
\leq d_\beta(\mathbf{F}_{j},0_{G^2(V)})^{4p} (Nt-Ns)^{4p\beta}
\]
for any $1/3<\beta<1/2$. Hypothesis \eqref{eq:bounddistancealpha} gives
\[
\Esp{\norm{X_{Ns}(s)^{-1}\cdot X_{Nt}(t)}^{4p}} \leq C (Nt-Ns)^{4p\beta} 
\]
for any $\beta<1/2$, hence the bound $N^{2p}(t-s)^{2p}$ by taking $\beta\to 1/2$.

If $\IntP{Ns}<\IntP{Nt}$, we introduce the event 
\begin{align*}
E = \left\{ T_{K(Ns)+1} < T_{K(Nt)} \right\}
\end{align*}
which corresponds to at least one complete excursion separating $Ns$ and $Nt$. On $E^c$, we use the decomposition
\[
X_{Ns}^{-1} \cdot X_{Nt} 
= 
\left(\mathbf{F}_{\IntP{Ns}}(\FracP{Ns})^{-1}
\cdot\mathbf{F}_{\IntP{Ns}}(1) \right)\cdot
\mathbf{F}_{\IntP{Ns}+1}(1)
\cdot\ldots\cdot
\mathbf{F}_{\IntP{Nt}-1}(1) \cdot
\mathbf{F}_{\IntP{Nt}}(\FracP{Nt})
\]
and, on $E$, we factorize the excursions using the $\ha{F}_k$:
\begin{subequations}
\begin{align}
\label{eq:beforeexc1}
X_{Ns}^{-1}\cdot X_{Nt} 
= &
\mathbf{F}_{\IntP{Ns}}(\FracP{Ns})^{-1}\cdot
\mathbf{F}_{\IntP{Ns}}(1)
\\
& \cdot \mathbf{F}_{\IntP{Ns}+1}(1)
\cdot\ldots\cdot
\mathbf{F}_{T_{K(Ns)+1}-1}(1)
\label{eq:beforeexc2}
\\
&\cdot \left(
\ha{F}_{K(Ns)+2}\cdot\ldots\cdot \ha{F}_{K(Nt))} 
\right)
\label{eq:duringexc}
\\
&
\cdot
\mathbf{F}_{T_{K(Nt)}}(1)
\cdot\ldots\cdot
\mathbf{F}_{\IntP{Nt}-1}(1)
\label{eq:afterexc1}
\\
& \cdot
\mathbf{F}_{\IntP{Nt}}(\FracP{Nt})
\label{eq:afterexc2}
\end{align}
\end{subequations}

\paragraph{First case: on $E$.}
We will use repeatedly the following ordering valid on $E$:
\[
T_{K(Ns)} \leq Ns < \IntP{Ns}+1 \leq T_{K(Ns)+1} < T_{K(Nt)} \leq \IntP{Nt} \leq Nt < T_{K(Nt)+1}
\]

We call $I_1$, $I_2$, $I_3$, $I_4$ and $I_5$ the respective norms of the terms \eqref{eq:beforeexc1}, \eqref{eq:beforeexc2},\eqref{eq:duringexc}, \eqref{eq:afterexc1} and \eqref{eq:afterexc2}. We then have
\[
\norm{\delta_{N^{-1/2}}(X_{Ns})^{-1} \cdot \delta_{N^{-1/2}}(X_{Nt})}^{4p}
\leq \frac{5^{4p}}{N^{2p}}  \sum_{k=1}^{5} I_k^{4p}
\]
Bounding $I_1^{4p}$ and $I_5^{4p}$ from above uses the distance~\eqref{eq:distancealpha}
\begin{align*}
\norm{ I_1}^{4p} & \leq  d_\beta(\mathbf{F}_{\IntP{Ns}},0_{G^2(V)})^{4p}(1-\FracP{Ns})^{4p\beta}
\\
\norm{ I_5}^{4p} & \leq  d_\beta(\mathbf{F}_{\IntP{Nt}},0_{G^2(V)})^{4p}(\FracP{Nt})^{4p\beta} 
\end{align*}
Taking $\sigma(R)$-conditional expectation with the bound~\eqref{eq:bounddistancealpha} gives a bound valid for all $\beta <1/2$:
\begin{align*}
\Espc{\norm{ I_1}^{4p}}{\sigma(R)} & \leq C (1-\FracP{Ns})^{4p\beta} \leq C (Nt-Ns)^{4p\beta}
\\
\Espc{\norm{ I_5}^{4p}}{\sigma(R)} & \leq C (\FracP{Nt})^{4p\beta} \leq C (Nt-Ns)^{4p\beta}
\end{align*}
and thus for $\beta \uparrow 1/2$ we obtain the expected bound $N^{2p}(t-s)^{2p}$ for $\Esp{\norm{I_1}^{4p}\indic{E}}$ and $\Esp{\norm{I_5}^{4p}\indic{E}}$.

The bounds on $I_2$ and $I_4$ are similar and we write down only the one for $I_2$ using again hypothesis~\eqref{eq:bounddistancealpha}:
\begin{align*}
\norm{I_2}^{4p} &\leq \left(\sum_{k=\IntP{Ns}+1}^{T_{K(Ns)+1}-1} d_\beta(\mathbf{F}_k,0_{G^2(V)})\right)^{4p}
\\
&\leq \left(T_{K(Ns)+1}-\IntP{Ns}-1\right)^{4p} \sum_{k=\IntP{Ns}+1}^{T_{K(Ns)+1}-1} d_\beta(\mathbf{F}_k,0_{G^2(V)})^{4p}
\\
&\leq C \left(T_{K(Ns)+1}-\IntP{Ns}-1\right)^{4p+1}
\\
&\leq C \left(T_{K(Ns)+1}-T_{K(Ns)}\right)^{2p+1} (Nt-Ns)^{2p}
\end{align*}
where the last inequality is true on $E$ only. The excursion times have finite moments from proposition~\ref{prop:inspectionmoments} and hypothesis~\ref{eq:finitemomentsT1} and thus:
\begin{equation}
\Esp{\norm{I_2}^{4p}\indic{E}} \leq C' (t-s)^{2p} 
\end{equation}

The bound on $I_3$ can be obtained using lemma~\ref{lemma:boundexcursionfortightness}. Since the number $K(Nt)-K(Ns)$ between $Ns$ and $Nt$ is necessarily smaller than $N(t-s)$, we have the bound
\begin{align*}
\Esp{ I_3(Ns,Nt)^{4p} \indic{E}}
\leq 
\Esp{ 
\sup_{1\leq k \leq N(t-s)} \norm{\ha{F}_{K(Ns)+2}\cdot\ldots\cdot\ha{F}_{K(Ns)+1+k}}^{4p}}
\end{align*}
We call $Z(T_{K(Ns)+1})$ the positive r.v. in the r.h.s. since it is
is a product of r.v. $\mathbf{F}_j(1)$ with $j\geq T_{K(Ns)+1}$. The filtration $(\ca{F}_n)$ is defined as $\ca{F}_n=\sigma((R_k,F_k);k\leq n)$. We now have:
\[
\Esp{Z(T_{K(Ns)+1}) } = \sum_{p\in\setN} \Esp{Z(T_{p+1})\indic{K(Ns)=p}}
=  \sum_{p\in\setN} \Esp{\Espc{Z(T_{p+1})}{\ca{F}_{T_{p+1}}}\indic{K(Ns)=p}}
\]
since the event $\{K(Ns)=p\} = \{T_p \leq Ns <T_{p+1}\}$ is $\ca{F}_{T_{p+1}}$-measurable. Using the strong Markov property for the hidden Markov chain $(R_n,\mathbf{F}_n)$ and the fact that $R_{T_{p+1}}=r_0$, we obtain
\[\Esp{Z(T_{K(Ns)+1}) }
= \sum_{p\in\setN} \Esp{ Z(0) \indic{K(Ns)=p}} = \Esp{Z(0)}
\]
We got rid of the dependency on $K(Ns)$ and lemma~\ref{lemma:boundexcursionfortightness} gives the final inequality
\[
\Esp{I_3(Ns,Nt)^{4p}\indic{E}} \leq B N^{2p} (t-s)^{2p}
\]

\paragraph{Second case: on $E^c$.}
On $E^c$, we have the following ordering:
\begin{equation}
\label{eq:orderingEc}
T_{K(Ns)} \leq \IntP{Ns} \leq Ns <\IntP{Ns}+1 \leq \IntP{Nt} \leq Nt < T_{K(Ns)+2}
\end{equation}
The difference is that the term $I_3$ is absent and $I_2$ and $I_4$ may be combined such that, on $E^c$:
\begin{align*}
\norm{X_{Ns}^{-1}\cdot X_{Nt} }^{4p} &\leq 3^{4p} \left(\norm{I_1}^{4p}+ \norm{I_{24}}^{4p} + \norm{I_5}^{4p}\right)
\end{align*}
where $I_{24}= \mathbf{F}_{\IntP{Ns}+1}(1)\cdot\ldots\cdot \mathbf{F}_{\IntP{Nt}-1}(1)$. The upper bounds on $I_1$ and $I_5$ are the same as in the previous case. For $I_{24}$, we have on $E^c$ 
\begin{align*}
\norm{ I_{24} }^{4p} \indic{E^c}
& \leq (\IntP{Nt}-\IntP{Ns}-1)^{4p} \sum_{k=\IntP{Ns}+1}^{\IntP{Nt}-1} 
d_\beta(\mathbf{F}_{k},0_{G^2(V)})^{4p}
\end{align*}
and thus now using \eqref{eq:bounddistancealpha} and the ordering~\eqref{eq:orderingEc} on $E^c$
\begin{align*}
\Espc{\norm{I_{24}}^{4p}}{\sigma(R)} \indic{E^c} 
& \leq C (\IntP{Nt}-\IntP{Ns}-1)^{4p+1} \indic{E^c}
\\
& \leq C (Nt-Ns)^{2p} (T_{K(Ns)+2}-T_{K(Ns)})^{2p+1} \indic{E^c}
\end{align*}
Hypothesis~\eqref{eq:finitemomentsT1} then gives the desired bound $C\Esp{T_2^{2p+1}}$ on $\Esp{\norm{I_{24}}^{4p}\indic{E^c}}$. 

Collecting all the previous bounds gives the expected tightness criterion since all the bounds are of the form $A N^{2p} (t-s)^{2p}$.
\end{proof}

\paragraph{Remarks on the hypotheses~\eqref{eq:finitemomentsT1} and \eqref{eq:bounddistancealpha}.}

Some of the hypotheses of theorem~\ref{thm:convergenceembeddings} could be slightly relaxed by improving the previous proof or by considering only a fixed given $\alpha<1/2$. In the case where only one value $\alpha<1/2$ is targeted, one could use the same approach as \citep{Breuillardetal} and require only a finite set of finite moments. However, in practice, one is often interested to the case $\alpha\to 1/2$. In this case, the finite moments hypothesis on $T_1$ cannot be relaxed since the variable $\ha{F}_k$ are required to have moments of all orders in order to apply the results of \citep{Breuillardetal}. Only efforts may be made on the requirement~\eqref{eq:bounddistancealpha} using for example correlations between the length $T_k$ of an excursion and the corresponding increment $\ha{F}_{k}$. The bound \eqref{eq:bounddistancealpha} is not restrictive in practice since it encompasses already the case where the embeddings $\mathbf{F}_k$ are smooth or Lipschitz.

\subsection{Iterated occupation times, the quasi-shuffle property and asymptotics}

\subsubsection{Proof of the quasi-shuffle property~\eqref{prop:occtime:quasishuffle}}
\label{sec:occtimeproof}
\begin{proof}
We first prove the following recursive decomposition. If $u=u_1\cdot u'$ with $u'=(u_2,\ldots,u_k)$ is the concatenation of the length one sequence $u_1$ and the sequence $u'$, then 
\begin{align*}
L_{u_1\cdot u'}(x) &= \sum_{0\leq i_1} \indic{x_{i_1}=u_1} \left( \sum_{i_1<i_2<\ldots<i_k\leq N} \indic{x_{i_2}=u_2}\ldots \indic{x_{i_k}=u_k} \right)
\\
&= \sum_{0\leq i_1 < N} \indic{x_{i_1}=u_1} L_{u'}( (x_{n})_{i_1+1\leq n< N})
\end{align*}
For two sequences $u=u_1u'$ and $v=v_1v'$ of length larger than $1$, the previous equation gives:
\begin{align*}
L_{u_1\cdot u'}(x)L_{v_1\cdot v'}(x) 
&= \sum_{0\leq i_1 < j_1 } \indic{x_{i_1}=u_1}\indic{x_{j_1}=v_1} L_{u'}((x_{n})_{i_1+1\leq n< N}) L_{v'}((x_{n})_{j_1+1\leq n< N})
\\
& \phantom{=} + \sum_{0\leq j_1 < i_1 } \indic{x_{i_1}=u_1}\indic{x_{j_1}=v_1} L_{u'}((x_{n})_{i_1+1\leq n< N}) L_{v'}((x_{n})_{j_1+1\leq n< N})
\\
&\phantom{=}+\sum_{0\leq i  } \indic{x_{i}=u_1}\indic{u_1=j_1} L_{u'}((x_{n})_{i+1\leq n< N}) L_{v'}((x_{n})_{i+1\leq n< N})
\\
&= \sum_{0\leq i_1 } \indic{x_{i_1}=u_1} L_{u'}((x_{n})_{i_1+1\leq n< N}) L_{v_1v'}((x_{n})_{i_1+1\leq n< N})
\\
& \phantom{=} + \sum_{0\leq j_1  } \indic{x_{j_1}=v_1} L_{u_1u'}((x_{n})_{j_1+1\leq n< N}) L_{v'}((x_{n})_{j_1+1\leq n< N})
\\
& \phantom{=} +\indic{u_1=v_1}\sum_{0\leq i  } \indic{x_{i}=u_1} L_{u'}((x_{n})_{i+1\leq n< N}) L_{v'}((x_{n})_{i+1\leq n< N})
\\
& = L_{u_1\cdot u'\qshuffle v_1\cdot v'}(x)
\end{align*}
from the definition of the quasi-shuffle product. The expected result is then obtained from the previous equation by recursion on the sum of the lengths of $u$ and $v$ by setting for convenience $L_\epsilon(x)=0$.
\end{proof}

\subsubsection{Relation between $C$, $\Gamma$ and the first iterated occupation times}
\label{sec:studyCGamma}

The covariance matrix $C$ and the anomalous area drift $\Gamma$ are obtained in property~\ref{prop:momentsofhatF} in terms of the moments of the law of the i.i.d.r.v. $\ha{F}_k$. It may be interesting to have more explicit formulae for $C$ and $\Gamma$, since they describe completely the limit law.

To this purpose, we introduce, for all $u\in E$, the following expectation values, which are related to the conditional law of the $\mathbf{F}_k$ and do not depend on the law of the process $(R_n)_{n\in\setN}$,
\begin{align*}
f_u &=\Espci{r_0}{\pi_1(\mathbf{F}_1(1)) }{R_1=u}
\\
c_u &=\Espci{r_0}{\pi_1(\mathbf{F}_1(1))\otimes \pi_1(\mathbf{F}_1(1))}{R_1=u} \in V\odot V
\\
\gamma_u  &= \Espci{r_0}{\pi_2(\mathbf{F}_1(1)) }{R_1=u} \in V\wedge V
\end{align*}
where $V\odot V$ is the symmetric subspace of $V\otimes V$.
The proof of proposition~\ref{prop:momentsofhatF} uses the fact that the invariant probability of the Markov chain $(R_n)$ satisfies, for any $u\in E$,
\[
\nu(u) = \frac{\Espi{r_0}{ L_u( (R_n)_{0\leq n<T_1} )}}{\Espi{r_0}{T_1}}
\]
Using the iterated occupation times, we now define, for any $k\in\setN$, the iterated measure on $E^k$:
\begin{equation}
\nu_{k,r_0}(u_1,\ldots,u_k) = \frac{\Espi{r_0}{ L_{u_1\ldots u_k}( (R_n)_{0\leq n<T_1} )}}{\Espi{r_0}{T_1}},
\end{equation}
which coincide, for $k=1$, with $\nu$. For $k\neq 1$, it depends on the initial point $r_0$: this can be seen for example in the total mass related to the moments of $T_1$, which depends on the initial point $r_0$. However, the expression below for $C$ and $\Gamma$ do not depend on $r_0$. It would be interesting to understand in more details the dependence on the initial point, but it is more a question of general theory of Markov processes than a rough path question.

\begin{property}
The covariance matrix $C$ and the anomalous area drift $\Gamma$, defined in property~\ref{prop:momentsofhatF} and appearing in theorems~\ref{thm:convergenceDonsker} and \ref{thm:convergenceembeddings} are given by:
\begin{subequations}
\begin{align}
\label{eq:formulaforCocctimes}
C_{ij} &= \sum_{u\in E} c_u^{(ij)} \nu(u) + \sum_{(u,v)\in E^2} (f_u^{(i)}f_v^{(j)}+f_u^{(j)}f_v^{(i)})\nu_2(u,v)
\\
\label{eq:formulaforGammaocctimes}
\Gamma_{ij} &= \sum _{u\in E} \gamma_u^{(ij)} \nu(u)
+ 
\frac{1}{2}\sum_{(u,v)\in E^2} (f_u^{(i)} f_v^{(j)}-f_u^{(j)} f_v^{(i)})\nu_2(u,v)
\end{align}
\end{subequations}
\end{property}
\begin{proof}
The proof is left to the reader and uses only the definition of a hidden Markov chain and the definition of $f_u$, $c_u$, $\gamma_u$, $\nu$ and $\nu_2$.
\end{proof}

The previous formula~\eqref{eq:formulaforGammaocctimes} for $\Gamma$ shows that there may be two ways of creating a non-zero $\Gamma$. The first way --- which is a bit trivial --- uses non-zero contributions $\gamma_u$ and corresponds to path increments $\mathbf{F}_k$ which already have a non-zero area $\pi_2(\mathbf{F}_k(1))$ in average: this is the case for example in the round-about model of section~\ref{subsec:examples}. However, this is impossible in the case of theorem~\ref{thm:convergenceDonsker} since geodesics in $V$ are straight lines with zero area.

The second way is much more interesting since it may create a non-zero anomalous area drift $\Gamma$ even in the context of theorem~\ref{thm:convergenceDonsker}: it is based on the $\nu_2$ contribution to $\Gamma$ in \eqref{eq:formulaforGammaocctimes}. In particular, it is absent from \citep{Breuillardetal}: random walks are a particular case of hidden Markov walk for which $E$ can be chosen to have cardinal $1$ and $T_1=1$ a.s. and thus $\nu_k=0$ for $k\geq 2$. One also checks easily that this term also vanishes for reversible Markov chain $(R_n)$, for which $\nu_2(u,v)=\nu_2(v,u)$ for all $(u,v)\in E^2$.

\subsubsection{Asymptotic fluctuations of the iterated occupation times and a remark on non-geometric rough paths}\label{sec:occtime:asymptotics}
Property~\ref{prop:ergodic:occtime} gives an almost sure limit of the rescaled quantities $L_u((R_n))_{n\in\setN}$; however, their fluctuations are more difficult to describe due to the quasi-shuffle property. We show here how to handle the question using our previous theorem~\ref{thm:convergenceembeddings}. To make things simpler, we focus on the case where $u$ is of length at most $2$ and $E$ is finite.

We now introduce the following quantities:
\begin{align*}
L^\star_{u_1}((R_n)_{0\leq n<N}) 
&= 
 \sum_{0\leq k <N} (\indic{R_k=u_1}-\nu(u_1))
 =L_{u_1}((R_n)_{0\leq n<N}) - N\nu(u_1) 
\end{align*}
If $E$ is finite, we introduce the finite-dimensional vector space $\setR^E$ with its canonical basis $(e_u)_{u\in E}$. We define the increments $F^\star_{k} = \sum_{u\in E}(\indic{R_k=u}-\nu(u)) e_u$ and the walk:
\begin{equation}
X^\star_N = \sum_{u\in E} L^\star_{u}((R_n)_{0\leq n<N})) e_u
\end{equation}

\begin{property}
The process $(R_n,X^\star_n)_{n\in\setN}$ is a hidden Markov walk in $E\times \setR^E$ with centred increments $(F^\star_n)$. If $(R_n)$ is irreducible, positive recurrent with finite moments of the first return time $T_1$, then theorem~\ref{thm:convergenceDonsker} may be applied to describe the scaling limit of the process.
\end{property}
This implies in particular, without any surprise, that $L_{u}((R_n)_{0\leq n<N})/N$ has (joint) Gaussian fluctuations of order $N^{-1/2}$ around its a.s. limit $\nu(u)$. But this also gives result about the second iterated occupation times $L_{u_1u_2}((R_n)_{0\leq n<N})$. To this purpose, we compute explicitly the iterated integral of the process $\iota_N(X^\star)$:
\begin{align*}
\frac{1}{N}\int_{0<s_1<s_2<1} d\iota_N(X^\star)(s_1)\otimes d\iota_N(X^\star)(s_2)
&= \frac{1}{2N}\sum_{k=0}^{N-1} F^\star_k\otimes F^\star_k \\
&
\phantom{=}+ \frac{1}{N} \sum_{0\leq k<l<N} F^\star_k\otimes F^\star_l
\end{align*}
The l.h.s. converges in law to $B_1\otimes B_1/2+A_{1}^{\text{L\'evy}} + \Gamma$. The first term of the r.h.s. converges a.s. and thus in law to a deterministic constant by the ergodic theorem.  The second term of the r.h.s. is related to the modified iterated occupation time through:
\begin{equation}
\label{eq:integralvssums}
\begin{split}
\frac{1}{N} \sum_{0\leq k <l < N} F_k^{(u_1)} F_l^{(u_2)} 
&= \frac{1}{N}\sum_{0\leq k<l<N} (\indic{R_k=u_1}-\nu(u_1))(\indic{R_l=u_2}-\nu(u_2))
\\
&=:
\frac{1}{N}
L^\star_{u_1u_2}((R_n)_{0\leq n <N}) 
\end{split}
\end{equation}
Relating $L^{\star}_{u_1u_2}$ to $L_{u_1u_2}$ gives the information about joint fluctuations of the collection of r.v. $L_{u_1u_2}((R_n)_{0\leq n<N})$ and additive functionals of the $L_u((R_n)_{0\leq n<N})$.

\paragraph{Remark on non-geometric rough paths.} The construction of the present section is purely combinatorial and involves only integer numbers for the indices: one may wonder why integrals of the Donsker embedding $\iota_N(X^\star)$ should be preferred to the choice of iterated sums for the signature in $T^2(V)$ (and not $G^2(V)$ any more) such as:
\[
\left( \sum_{0\leq k<N}F^\star_k, \sum_{0\leq k < l <N} F^\star_k\otimes F^\star_l\right)
\]
There is indeed no reason to prefer one to another and it may depend on the context. The positive result is that there is no reason since both constructions differ by a term $(1/2)\sum_{0\leq k<N} F^\star_k\otimes F^\star_k$, which belongs to the center of $T^2(V)$, is symmetric and whose limit is governed in the present case by the law of large numbers and is given by $(0,K t)$ where $K$ is a deterministic symmetric matrix. This additional term breaks the geometric rough path property but, as emphasized in exercise 2.14 (page 23) of \citep{FrizHairerBook} (and remarks disseminated in the corresponding chapter) or \citep{artHaiKelly}, this subtlety does not make a big difference from an analytic perspective.

\section{Some open questions and extensions}

We have seen how to build non-trivial rough paths above Brownian motion from very simple and intuitive processes such as hidden Markov chains. This raises various questions.

A key role in the emergence of a non-zero area anomaly $\Gamma$ is played by the short-time correlations of the underlying Markov chain $(R_n)_{n\in\setN}$. Exact renormalization on the time scale is due to the excursion decomposition: it would be interesting to generalize it to more general processes, such as $\alpha$-mixing processes as described in \citep{bookBillConv}.

On one hand, we have put restrictive hypothesis on the return times $T_1$ and the moments of the increments $\mathbf{F}_n$, so that the limit belongs to the Brownian universality class. On the other hand, rough paths structure may also describe L\'evy processes as described in \citep{ChevyrevLevy}. It would be interesting to build discrete time models that converge to such L\'evy processes and contain all types of admissible anomalies such as $\Gamma$.

The generalization of theorem~\ref{thm:convergenceDonsker} to theorem~\ref{thm:convergenceembeddings} uses embeddings. This question of a discrete structure on top of piecewise paths is similar to the theory of piecewise-deterministic Markov processes (PDMP) in continuous time and it may be interesting to study space-time renormalized PDMP using the present rough path approach.

Rough paths are a particular case of much more general regularity structures as introduced by \citep{FrizHairerBook}. One may expect that such regularity structures may contain a wide class of anomalies (both in the sense of our area anomaly and in the sense of anomalies in field theory, i.e. a broken symmetry in the discretization or regularization restored by counter-terms in the continuous limit) and it may interesting to understand them from correlations in discrete models, as done in the present paper. In particular, anomalous enhanced Brownian motion is a particular case of translation of rough paths as described in \citep{Roughpathrenormalization}: the renormalization scheme described in this reference corresponds to our excursion-based renormalization. It would also be interesting in this context to introduce branched version of our iterated occupation times.

The other novelty is the introduction of iterated occupation time and the emergence of shuffle or quasi-shuffle products already at the discrete level. Such products also appear in other domains of algebra or combinatorics, for example in the theory of multiple zeta values or periods: it would be interesting to examine whether relevant Markov chains could be related to such theories.

\section*{Acknowledgements}
D.~S. is partially funded by the Grant ANR-14CE25-0014 (ANR GRAAL). D.S. and O.L. thank Lorenzo Zambotti for stimulating discussions and useful suggestions.

\bibliographystyle{imsart-nameyear}

\end{document}